\documentclass[12pt]{amsart}
\usepackage{amssymb}
\newtheorem{theorem}{Theorem}
\newtheorem{proposition}[theorem]{Proposition}

\newtheorem{lemma}[theorem]{Lemma}
\newtheorem{example}[theorem]{Example}

\title{Upper bound for the Lempert function of smooth domains}

\author[N.~Nikolov, P.~Pflug, P.~J.~Thomas]
{Nikolai Nikolov, Peter Pflug, Pascal J.~Thomas}

\address
{Institute of Mathematics and Informatics\\ Bulgarian Academy of
Sciences\\ Acad. G. Bonchev 8, 1113 Sofia,
Bulgaria}\email{nik@math.bas.bg}

\address{Carl von Ossietzky Universit\"at Oldenburg\\
Institut f\"ur Mathe\-ma\-tik\\ Postfach 2503\\ D-26111 Oldenburg,
Germany}\email{pflug@mathematik.uni-oldenburg.de}

\address{Institut de Math\'ematiques\\
Universit\'e Paul Sabatier, 118 Route de Narbonne, 31062 Toulouse
Cedex 9, France} \email{pthomas@cict.fr}



\begin{document}

\begin{thanks}{This paper was written during the stay of the first-named
author at the Universit\'e Paul Sabatier, Toulouse (October 2008)
supported by a CNRS programme "Convention d'\'echanges" No 21463
and the stay of the second-named author at the Institute of
Mathematics and Informatics of the Bulgarian Academy of Sciences
(September 2008) supported by a DFG grant 436POL113/103/0-2.}
\end{thanks}

\begin{abstract} An upper estimate for the Lempert function of any
$C^{1+\varepsilon}$-smooth bounded domain in $\Bbb C^n$ is found
in terms of the boundary distance.
\end{abstract}

\maketitle

\section{Introduction}

By $\Bbb D$ we denote  the unit disc in $\Bbb C.$ Let $D$ be a
domain in $\Bbb C^n.$ Recall the definition of the Lempert
function of $D:$
$$
l_D(z,w)=\inf\{\alpha\in[0,1):\exists\varphi\in\mathcal O(\Bbb
D,D):\varphi(0)=z,\varphi(\alpha)=w\},
$$
where $z,w\in D.$ The Kobayashi pseudodistance $k_D$ is the
largest pseudodistance below $\tilde k_D=\tanh^{-1}l_D$.

Combining the proofs of Proposition 2.5 in \cite{FR} and
Proposition 10.2.3 in \cite{JP}, it follows that if $D$ is bounded
and $C^{1+\varepsilon}$-smooth, then there exists $c>0$ such that
$$(\ast)\quad
2k_D(z,w)\le\log\left(1+\frac{||z-w||}{d_D(z)}\right)+\log
\left(1+\frac{||z-w||}{d_D(w)}\right)+c,$$ where $d_D$ is the
distance to $\partial D.$ In \cite{FR} this inequality is applied
for extension of proper holomorphic maps.

We point out that the assumption of smoothness is essential; the
conclusion fails if $D$ is a planar polygon (its boundary is
Lipschitz).

On the other hand, one can show that if $D$ is strongly
pseudoconvex, then for any $\delta>0$ there exists a $c'>0$ such
that
$$2k_D(z,w)\ge 2c_D(z,w)\ge-\log d_D(z)-\log d_D(w)-c'\quad\hbox{if
}||z-w||\ge\delta,$$ where
$c_D(z,w)=\sup\{\tanh^{-1}|f(w)|:f\in\mathcal O(D,\Bbb
D),f(z)=0\}$ is the Ca\-rath\'eodory distance of $D.$ It is
natural to ask whether the reverse inequality holds for $\tilde
k_D$ (see \cite{NPT}). Then the reverse inequality will imply
$(\ast)$ if $z$ and $w$ are not close to one and the same boundary
point. Such an estimate is equivalent to the following one.

\begin{theorem}\label{1} Let $D\subset\Bbb C^n$ be a $C^{1+\varepsilon}$-smooth
bounded domain. Then there is $c>0$ such that for any $z,w\in D,$
$$l_D(z,w)\le 1-cd_D(z)d_D(w).$$
\end{theorem}

The following example may show that the smoothness assumption in
the previous theorem is important.

\begin{example}\label{4} Let $D\subset\Bbb C$ be the image of $\Bbb D$
under the map $z\to 2z+(1-z)\log(1-z).$ Then $D$ is a $C^1$-smooth
bounded domain and
$$\lim_{\Bbb R\ni w\uparrow2}\frac{1-l_D(0,w)}{d_D(w)}=0.$$
\end{example}

The main difficulty in the proof of Theorem \ref{1} arises from
the fact that, in general, $l_D$ does not satisfy the triangle
inequality. Therefore, we cannot localize this function,
i.e.~reduce the proof to the case when both arguments are near the
boundary. On the other hand, the idea for the proof is clear: join
two points (possibly on the boundary) by a suitable real-analytic
curve in the domain and perturb a holomorphic extension of the
curve in order to cover neighborhoods of these points. To control
the perturbation, we shall need the upper bound in the following
proposition.

\begin{proposition}\label{2} Let $(G_j)$ be a
sequence of $C^1$-smooth simply connected bounded planar domains
which tend to a bounded domain $G\subset\Bbb C$ in a sense that
for any two open sets $K,L\subset\Bbb C$ with $K\Subset G\Subset
L$ there exists an $j_0$ such that $K\subset G_j\subset L$ for any
$j>j_0.$ Let $\varepsilon>0.$ Assume that the following regularity
properties hold: there is a neighborhood $U$ of $\partial G$ such
that
$$||\nabla r_j(z)||\ge\varepsilon,\ ||\nabla r_j(z)-\nabla
r_j(w)||\le\varepsilon^{-1}||z-w||^{\varepsilon} ,\quad j\in\Bbb
N,\ z,w\in U,$$ where $r_j\in C^1(\Bbb C)$ is a defining function
of $G_j$ {\rm (}i.e. $G_j=\{z\in\Bbb C: r_j(z)<0\}${\rm )}. Then
there is a $c>0$ such that for any conformal map $f_j:\Bbb D\to
G_j$ with $d_{G_j}(f_j(0))\ge\varepsilon$ one has that
$$c^{-1}\le|f'_j(z)|\le c,\quad z\in\Bbb D.$$
\end{proposition}

\noindent{\bf Remark.} Recall that if $D$ is the inner domain of a
Dini-smooth closed Jordan curve $\gamma$ \footnote{This means that
$\int_0^1\frac{\omega(t)}{t}dt<\infty,$ where $\omega$ is the
modulus of continuity of $\gamma.$}, then any conformal map
$f:\Bbb D\to D$ extends to a diffeomorphism from $\overline{\Bbb
D}$ to $\overline D$ (cf. Theorem 3.5 in \cite{Pom}). Since
$\tilde k_D=k_D$ for any $D\subset\Bbb C,$ it is not difficult to
see that Theorem 1 holds for $n=1,$ if $\partial D$ is a
Dini-smooth curve near any boundary point.

\section{Proofs}

\noindent{\it Proof of Example \ref{4}.} The facts that the map is
injective on $\Bbb D$ and $G$ is a $C^1$-smooth domain may be
found in \cite{Pom}, p.~46 (use e.g. Proposition 1.10). Let
$\psi:D\to\Bbb D$ be the inverse map. Then
$$l_D(0,w)=l_{\Bbb D}(0,\psi(w))=|\psi(w)|.$$
On the other hand, since $D$ is $C^1$-smooth, it is not difficult
to see that (use e.g. Proposition 2 in \cite{JN})
$$\lim_{w\to\partial D}|\psi'(w)|\frac{d_D(w)}{d_{\Bbb
D}(\psi(w))}=1.$$ Hence,
$$\lim_{\Bbb R\ni w\uparrow2}\frac{1-l_D(0,w)}{d_D(w)}=
\lim_{\Bbb R\ni w\uparrow2}|\psi'(w)|=0.$$

\noindent{\it Proof of Proposition \ref{2}.} All the constant
below will be independent of $j$ and of the boundary points that
appear.

Let $D\subset\Bbb C$ be a $C^1$-smooth bounded domain. For any
point $a\in\partial D$ there is a $\theta_a\in\Bbb R$ such that if
$\rho_a:z\to(z-a)e^{i\theta_a}$ and $D_a=\rho_a(D),$ then $x>0$
($z=x+iy$) is the inward normal to $\partial D_a$ at $0.$
Moreover, for $\delta>0$ put $G^i_\delta=\{z\in\Bbb C:
|z|<2\delta, x>|y|^{1+\delta}\}$ and $G^e_\delta=\Bbb
C\setminus\overline{(-G^i_\delta)}.$

Using the assumptions of the proposition, we may shrink
$\varepsilon$ such that for any $j$ and any $a\in\partial G_j$ one
has that
$$G^i_\varepsilon=:G^i\subset G_{j,a}\subset G^e:=G^e_\varepsilon.$$
We may smooth $G^i$ and $G^e$ at their  angular points preserving
these inclusions.

To prove the result, we shall need the following two estimates.
\smallskip

\noindent{\it Estimate 1.} There exists a $c_1>0$ such that
$$
c_1^{-1}\le \kappa_{G_j}(z;1)d_{G_j}(z)\le c_1,\quad z\in
G_j,\;j\in\Bbb N,
$$
where $\kappa_D$ denotes the Kobayashi-Royden metric of an
arbitrary domain $D\subset\Bbb C. $\footnote{Recall that
$\kappa_D(z;X)=\inf\{\alpha\ge 0:\exists\varphi\in\mathcal O(\Bbb
D,D):\varphi(0)=z,\alpha\varphi'(0)=X\}$ for any domain
$D\subset\Bbb C^n.$}

\noindent{\it Subproof.} Fix $j\in\Bbb N$ and $z\in G_j$. First,
we shall prove the lower bound. Let $\psi$ be a conformal map from
$G^e$ to $\Bbb D_\ast.$ Choose $a\in\partial G_j$ such that
$d_{G_j}(z)=|z-a|$ and put $z_a=\rho_a(z).$ Then
$$\kappa_{G_j}(z;1)\ge\kappa_{G^e}(z_a;1)=\kappa_{\Bbb D_\ast}(\psi(z_a);\psi'(z_a))
\ge\kappa_{\Bbb D}(\psi(z_a);\psi'(z_a))$$ $$\ge c_3 d^{-1}_{\Bbb
D}(\psi(z_a))\ge c_4d_{G_e}^{-1}(z_a)=c_4d^{-1}_{G_j}(z)$$ (the
constant $c_3$ is provided by the facts that $\psi$ extends to a
diffeomorphism from $\overline{G^e}$ to $\overline{\Bbb D}_\ast$
and that $\cup_{j=1}^\infty G_j$ is a bounded set).

Next, we prove the upper bound. Assume the contrary. Then we may
find a sequence of points $z_j\in G_j$ (if necessary take an
appropriate subsequence) such that
$\kappa_{G_j}(z_j;1)d_{G_j}(z_j)\to\infty.$ Since for any domains
$G^2\Subset G^1\Subset G$ one has that
$$\kappa_{G_j}(z;1)\le\kappa_{G^1}(z;1)\le c(G^2),\quad j\gg1,
z\in G^2,$$ it follows that $d_{G_j}(z_j)\to 0.$ Let
$a_j\in\partial G_j$ be such that $d_{G_j}(z_j)=|z_j-a_j|$ and
$z_{a_j}=\rho_{a_j}(z_j).$ Then $z_{a_j}\in
G_i\cap(0,\varepsilon)$ for any $j\gg 1.$ Similar to above, we get
that
$$\kappa_{G_j}(z_j;1)\le\kappa_{G^i}(z_{a_j};1)\le c_5d_{G^i}^{-1}(z_{a_j})
\le c_6z^{-1}_{a_j}=c_6d^{-1}_{G_j}(z_j),$$ a contradiction.
\smallskip

\noindent{\it Estimate 2.} There is a $c_2>0$ such that if
$d_{G_j}(w)\ge\varepsilon,$ then
$$c_2^{-1}\le\frac{1-l_{G_j}(z,w)}{d_{G_j}(z)}\le c_2.$$

\noindent{\it Subproof.} Following the proof of the lower bound in
Estimate 1, we have
$$l_{G_j}(z,w)\ge l_{G^e}(z_a,w_a)=l_{\Bbb D}(\psi(z_a),\psi(w_a))$$
$$\ge 1-c_7d_{\Bbb D}(\psi(z_a))\ge 1-c_8d_{G^e}(z_a)=1-c_8d_{G_j}(z)$$
and the upper bound is proved (not using that
$d_{G_j}(w)\ge\varepsilon$).

Next, we shall prove the lower bound. Since
$d_{G_j}(w)\ge\varepsilon,$ following the proof of the upper bound
in Estimate 1, it is enough to show the bound when $z_a\in
G^i\cap(0,\varepsilon).$ Then
$$k_{G_j}(z,w)=k_{G_{j,a}}(z_a,w_a)\le k_{G_{j,a}}(z_a,\varepsilon)+
k_{G_{j,a}}(\varepsilon,w_a)$$
$$\le k_{G_{j,a}}(z_a,\varepsilon)+c_9
\le k_{G^i}(z_a,\varepsilon)+c_9$$ (to see the second inequality,
use $G$ as above) and hence
$$1-l_{G_j}(z,w)\ge c_{10}(1-l_{G^i}(z_a,\varepsilon))\ge c_{11}
d_{G^i}(z_a)\ge c_{12}d_{G_j}(z).$$

Now, using both estimates, we shall prove the desired
inequalities. Since $|z|=l_{\Bbb D}(0,z)=l_{G_j}(f_j(0),f_j(z)),$
it follows by Estimate 2 that
$$c_2^{-1}\le\frac{d_{G_j}(f_j(z))}{d_{\Bbb D}(z)}\le c_2.$$
On the other hand,
$$(1-|z|^2)^{-1}=\kappa_{\Bbb
D}(z;1)=\kappa_{G_j}(f_j(z),f_j'(z))=|f_j'(z)|\kappa_{G_j}(f_j(z);1).$$
Then Estimate 1 implies that
$$c_1^{-1}\le\frac{d_{G_j}(f_j(z))}{|f_j'(z)|d_{\Bbb D}(z)}\le 2c_1.$$
Hence $(2c_1c_2)^{-1}\le|f_j'(z)|\le c_1c_2.$
\medskip

\noindent{\it Proof of Theorem \ref{1}.} By compactness, it is
enough to prove the estimate when $z$ and $w$ are near two
boundary points $a$ and $b$ (possible $a=b$), and when $z$ lies in
a compact subset of $D,$ but $w$ is near a boundary point $b.$ We
shall consider only the first case, because the second one is
similar and even simpler.

\begin{lemma}\label{3} There is  a
polynomial map $\varphi:\Bbb C\to\Bbb C^n$ such that
$$\varphi((-1,1))\subset D,\varphi(1)=a,\varphi(-1)=b,
\varphi'(1)=-n_a,\varphi'(-1)=n_b,$$ where $n_p$ is the inward
normal vector to $\partial D$ at $p.$
\end{lemma}

For completeness we shall prove this lemma at the end the paper.
\smallskip

Let $(u,v)\in T_a^{\Bbb C}\partial D\times T_b^{\Bbb C}\partial
D.$ Set
$$\varphi_{u,v}(\zeta)=\varphi(\zeta)+\left(\frac{\zeta+1}{2}\right)^2u+
\left(\frac{\zeta-1}{2}\right)^2v, \quad\zeta\in\Bbb C;$$
$$\Phi(\zeta_1,u,\zeta_2,v)=(\varphi_{u,v}(\zeta_1),\varphi_{u,v}(\zeta_2)).$$
Computing the Jacobian, it follows that $\Phi:\Bbb C^{2n}\to\Bbb
C^{2n}$ is invertible at $(1,0,-1,0).$ Then for $(z,w)$ in a
neighborhood $U\subset\Bbb C^{2n}$ of $(a,b)$ there are
$\zeta_1(z,w),\zeta_2(z,w),u(z,w),v(z,w)$ such that
$$z=\varphi_{u(z,w),v(z,w)}(\zeta_1(z,w)),\
w=\varphi_{u(z,w),v(z,w)}(\zeta_2(z,w)).$$ Set
$\psi_{z,w}=\varphi_{u(z,w),v(z,w)}$ and
$G_{z,w}=\psi_{z,w}^{-1}(D).$ If $r$ is a
$C^{1+\varepsilon}$-smooth defining function of $D,$ then
$\rho_{z,w}=r\circ\psi_{z,w}$ is a defining function of $G_{z,w}.$
Note that $(1,0)$ and $(-1,0)$ are the outward normal vectors to
$\partial G_{a,b}$ at the points $A=(1,0)$ and $B=(-1,0),$
respectively. Then we may shrink $U$ and find two squares
$S_\delta(A)$ and $S_\delta(B)$ with side $2\delta$ and centers at
$A$ and $B,$ respectively, such that $\rho_{z,w}$ is close (in
sense of the $C^{1+\varepsilon}$-norm) to the function $x-1$ in
$S_\delta(A)$ and  to the function $-x-1$ in $S_\delta(B).$ Since
$[-1+\delta/2, 1-\delta/2]\Subset G_{a,b}$, we may find
$\delta'>0$ such that
$R_{-\delta,\delta'}=\{\zeta:|x|<1-\delta/2,|y|<\delta'\}\Subset
G_{a,b}.$ We may shrink $U$ such that $R_{-\delta, \delta'}\Subset
G_{z,w}$ for any $(z,w)\in U.$ Shrinking $U$ further, we may
assume that the curve $\gamma_{z,w}=\{\rho_{z,w}=0\}$ intersects
only once each of the horizontal line segments of length $2\delta$
inside $S_\delta (A)$, and likewise for $S_\delta (B)$. Then
$H_{z,w}=G_{z,w}\cap R_{2\delta,\delta'}$ is bounded by two
horizontal line segments contained in $\{|y|= \delta'\}$ and by
the curves $\gamma_{z,w}\cap S_\delta (A) $ and $\gamma_{z,w}\cap
S_\delta (B);$ so $H_{z,w}$ fails to be smooth only at its four
corners. We smooth $H_{z,w}$ such that it remains unchanged
outside of a $\delta'/2$ neighborhood of the corners, and that
$H_{z,w}$ is close to $H_{a,b}$ (as before). Let $\eta_{z,w}:\Bbb
D\to H_{z,w}$ be a conformal map with $\eta_{z,w}(0)=0$ and
$\eta_{z,w}(p_{j,z,w})=\zeta_j(z,w),$ $j=1,2.$ It extends to a
diffeomorphism from $\overline{\Bbb D}$ to $\overline {H_{z,w}}.$
By Proposition \ref{2}, reasoning by contradiction, we may shrink
$U$ and $\varepsilon$ such that $\eta'_{z,w}$ are uniformly
bounded from above. Setting $q_{j,z,w}=p_{j,z,w}/|p_{j,z,w}|$ and
shrinking $U$ once more, it follows by the mean-value inequality
that $\rho_{z,w}(q_{j,z,w})\in\partial G_{z,w}.$ The same
inequality for $\theta_{z,w}=\psi_{z,w}\circ\eta_{z,w}\in\mathcal
O(\Bbb D,D)$ implies that
$$d_D(z)\le||\theta_{z,w}(q_{1,z,w})-\theta_{z,w}(p_{1,z,w})||\le
Cd_{\Bbb D}(p_{1,z,w})$$ (since
$\theta_{z,w}(q_{1,z,w})\in\partial D$) and similarly $d_D(w)\le
Cd_{\Bbb D}(p_{2,z,w}).$ Hence
$$1-l_D(z,w)\ge1-\left|\frac{p_{1,z,w}-p_{2,z,w}}{1-p_{1,z,w}
\overline{p_{2,z,w}}}\right|>\frac{d_{\Bbb D}(p_{j,z,w})d_{\Bbb
D}(p_{j,z,w})}{2}\ge\frac{d_D(z)d_D(w)}{2C^2}.$$
\medskip

\noindent{\it Proof of Lemma \ref{3}.} In the proof we will only
assume that $D$ is $C^1$-smooth near $a$ and $b.$ We start with a
$C^2$-smooth curve $\tilde \varphi:[-1,1]\to\Bbb C^n$ such that
$$\tilde\varphi((-1,1))\subset D,\tilde
\varphi(1)=a,\tilde\varphi(-1)=b,
\tilde\varphi'(1)=-n_a,\tilde\varphi'(-1)=n_b.$$ Then for
$\varepsilon>0$ choose a polynomial map
$\varphi_{\varepsilon}:\Bbb C\to\Bbb C^n$ that agrees with
$\tilde\varphi$ at $\pm 1$ up to order 1 and such that
$||\varphi'_{\varepsilon}(t)-\varphi'(t)||<\varepsilon$ for any
$t\in(-1,1)$. This map will do the job for any small
$\varepsilon.$

Indeed, we shall show that there are $\varepsilon_1,\delta_1>0$
such that $\varphi_\varepsilon((1-\delta_1,1))\subset D$ for any
$\varepsilon<\varepsilon_1.$ Let $r$ be a defining function of $D$
which is $C^1$-smooth near $a$ and $b.$ Put
$\rho_\varepsilon=r\circ\varphi_\varepsilon.$ Then there exists a
$\delta\in(0,1)$ such that
$$\rho_\varepsilon(1-t)=-2\int_{1-t}^1\hbox{Re}\langle\partial
r(\varphi_\varepsilon(s)),\overline{\varphi'_\varepsilon(s)}\rangle
ds,\quad 0<t<\delta.$$ Shrinking $\delta,$ we may assume that
$||2\partial r(z)+n_a||<1/4$ if $||z-a||<\delta;$ in particular,
$||\partial r(z)||<5/8.$ Since
$$\varphi_\varepsilon(1-t)=a-\int_{1-t}^1\varphi'_\varepsilon(s)ds$$
and $||\varphi'_\varepsilon(s)-\tilde\varphi'(s)||<\varepsilon,$
there are $\varepsilon_1,\delta_1>0$ such that
$||\varphi_\varepsilon(s)-a||<\delta$ and
$||\varphi'_\varepsilon(s)+n_a||<1/5$ if $1-s<\delta_1$ and
$\varepsilon<\varepsilon_1.$ Thus,
$$|1-2\hbox{Re}\langle\partial
r(\varphi_\varepsilon(s)),\overline{\varphi'_\varepsilon(s)}\rangle|\le
2|\hbox{Re}\langle\partial
r(\varphi_\varepsilon(s)),\overline{\varphi'_\varepsilon(s)+n_a}\rangle|$$
$$+|\hbox{Re}\langle2\partial
r(\varphi_\varepsilon(s))+n_a,\overline{n_a}\rangle|<\frac{5}{4}.\frac{1}{5}+
\frac{1}{4}.$$ Hence, $\hbox{Re}\langle
r(\varphi_\varepsilon(s)),\overline{\varphi'_\varepsilon(s)}\rangle
>1/4,$ which implies that
$\rho_\varepsilon(1-t)<-t/2,$ and we are done.

Similarly,  there exist $\varepsilon_2,\delta_2>0$ such that
$\varphi_\varepsilon((-1,-1+\delta_2))\subset D$ for any
$\varepsilon<\varepsilon_2.$ Note that for
$\delta_3=\min\{\delta_1,\delta_2\}$ there is an $\varepsilon_3>0$
such that $\varphi_\varepsilon([-1+\delta_3,1-\delta_3])\subset D$
for any $\varepsilon<\varepsilon_3.$ Therefore, any
$\varepsilon<\min\{\varepsilon_1,\varepsilon_2,\varepsilon_3\}$
does the job.

\end{document}